\documentclass[12pt]{article}

\usepackage{amsfonts}
\usepackage{amssymb,amsmath}
\usepackage{amsthm}

\usepackage{todonotes}

\usepackage[upright]{fourier}\usepackage{baskervald}
\usepackage{dsfont}

\usepackage{graphicx}
\usepackage{float}
\usepackage{wrapfig}
\usepackage{subfigure}

\usepackage{url}

\usepackage{hyperref}
\usepackage{geometry}

\usepackage{verbatim}

\usepackage[english]{babel}
\usepackage[latin1]{inputenc}

\usepackage{calligra}
\usepackage[T1]{fontenc}

\usepackage{bbm}

\usepackage{mathrsfs}

\usepackage{color}
\usepackage[color,all]{xy}

\theoremstyle{definition}
\newtheorem{definition}{Definition}[section]
\newtheorem{example}[definition]{Example}
\newtheorem{remark}[definition]{Remark}

\theoremstyle{plain}
\newtheorem{theorem}[definition]{Theorem}

\newtheorem{conjecture}{Conjecture}

\newtheorem{proposition}[definition]{Proposition}

\theoremstyle{remark}

\newcommand{\OO}{\mathscr{O}}

\newcommand{\PP}{\mathbb{P}}

\newcommand{\codim}{\text{\upshape{codim}}}

\newcommand{\Pic}{\text{\upshape{Pic}}}

\newcommand{\cE}{\mathscr{E}}

\newcommand{\PPE}{{\PP(\cE)}}

\newcommand{\cL}{\mathscr{L}}
\newcommand{\cH}{\mathscr{H}}

\newcommand{\cM}{\mathscr{M}}
\newcommand{\cHom}{\mathscr{H}\!om}
\newcommand{\rk}{\text{\upshape{rk}}}
\newcommand{\QQ}{\mathbb{Q}}
\newcommand{\ZZ}{\mathbb{Z}}

\makeatletter
\newcommand{\subjclass}[2][2010]{%
  \let\@oldtitle\@title%
  \gdef\@title{\@oldtitle\footnotetext{#1 \emph{Mathematics subject classification:} #2}}%
}
\newcommand{\keywords}[1]{%
  \let\@@oldtitle\@title%
  \gdef\@title{\@@oldtitle\footnotetext{\emph{Key words and phrases:} #1.}}%
}
\makeatother

\AtEndDocument{\bigskip{
  \textsc{Universit\"at des Saarlandes, Campus E2 4, D-66123 Saarbr\"ucken, Germany} \par  
  \textit{E-mail address}, C.~Bopp: \texttt{bopp@math.uni-sb.de} \par
  \textit{E-mail address}, M.~Hoff: \texttt{hahn@math.uni-sb.de} \par
}}

\begin{document}

\title{The relative canonical resolution: \texttt{Macaulay2}-package, experiments and conjectures}
\author{Christian Bopp and Michael Hoff}
\date{}


\maketitle

\begin{abstract}
This short note provides a quick introduction to relative canonical resolutions of curves on rational normal scrolls. We present our \texttt{Macaulay2}-package which computes the relative canonical resolution associated to a curve and a pencil of divisors. Most of our experimental data can be found on the following webpage \cite{M2exp}. We end with a list of conjectural shapes of relative canonical resolutions. In particular, for curves of genus $g=n\cdot k +1$ and pencils of degree $k$ for $n\ge 1$, we conjecture that the syzygy divisors on the Hurwitz scheme $\cH_{g,k}$ constructed in \cite{DP18} all have the same support. 
\end{abstract}

\section{Relative Canonical Resolutions} \label{relCanRes}

The \emph{relative canonical resolution} is the minimal free resolution of a canonically embedded curve $C$ inside a rational normal scroll. Every such scroll is swept out by linear spaces parametized by pencils of divisors on $C$.  

Studying divisors on moduli spaces reveals certain aspects of the global geometry of these spaces. 
A famous example for curves of odd genus $g$ is the \emph{Koszul-divisor} on the moduli space of curves $\cM_g$. It is induced by the minimal free resolution of $C\subset \PP^{g-1}$. Set-theoretically the Koszul-divisor consists of curves, such that the minimal free resolution of the canonical model has extra-syzygies at a certain step. 
In \cite{BP},\cite{DP} and \cite{DP18}, the relative canonical resolution was used to define similar \emph{syzygy divisors} on Hurwitz spaces $\cH_{g,k}$, parametizing pairs of curves of genus $g$ and pencils of divisors of degree $k$  (equivalently, covers of $\PP^1$ of degree $k$ by curves of genus $g$).    

We will briefly summarize the connections between pencils of divisors on canonical curves and 
rational normal scrolls in order to define the relative canonical resolution. 
Let $C \subset \PP^{g-1}$ be a canonically embedded curve of genus $g$ and let further 
$$g^1_k=\{D_\lambda\}_{\lambda \in \PP^1}\subset |D|$$ be a pencil of divisors of degree $k$. 
If we denote by $\overline{D_\lambda}\subset \PP^{g-1}$ the linear span of the divisor, then
$$
X=\bigcup_{\lambda\in \PP^1}\overline{D_\lambda}\subset \PP^{g-1}
$$
is a $(k-1)$-dimensional rational normal scroll of degree $f=g-k+1$. 

\begin{definition}
Let $e_1 \geq e_2 \geq \dots \geq e_d\geq 0$ be integers, 
$\cE=\OO_{\PP^1}(e_1)\oplus\dots\oplus\OO_{\PP^1}(e_d)$ and let 
$\pi: \PPE\to \PP^1$ be the corresponding $\PP^{d-1}$-bundle.\newline
A \emph{rational normal scroll} $X=S(e_1,\dots,e_d)$ of  type $(e_1,\dots,e_d)$ is the image of 
$$
j:\PP(\cE)\to \PP H^0(\PPE,\OO_\PPE(1))= \PP^r
$$ \vspace{-1mm}
where $r=f+d-1$ with $f=e_1+\dots+e_d\geq 2$.
\end{definition}

Conversely if $X$ is a rational normal scroll of degree $f$ containing a canonical curve, 
then the ruling on $X$ cuts out a pencil of divisors $\{D_\lambda\}\subset |D|$ such that $h^0(C,\omega_C\otimes \OO_C(D)^{-1})=f.$

\begin{example}
We consider a non-hyperelliptic canonically embedded curve $C\subset \PP^3$ of genus $4$. The curve $C$ is a complete intersection of a quadric surface $Q$ and a cubic surface $S$. If $C$ admits exactly two pencils of degree $3$ (which is also the maximal number), then the quadric $Q$ is isomorphic to $\PP^1\times \PP^1$. By B\'ezouts theorem, the two rulings of lines on $Q$ cut out the two pencils of degree $3$ on $C$ and conversely, the quadric is the scroll of type $(1,1)$ swept out by any of these pencils. If $C$ only admits one pencil of degree $3$, then   the quadric $Q$ is isomorphic to a cone (i.e., a quadric of rank $3$) and coincides with the scroll of type $(0,2)$ swept out by the unique pencil.
\end{example}

In \cite{H} it is shown that the variety $X$ defined above is a non-degenerate $d$-dimensional variety of minimal degree
$\deg X=f=r-d+1=\codim X+1$. If  $e_1,\dots,e_d>0$, then $j:\PP(\cE)\to X\subset \PP H^0(\PPE,\OO_\PPE(1))= \PP^r$ is an isomorphism. Otherwise, it is a resolution of singularities. Since $R^i j_* \OO_\PPE=0$, it is convenient to consider $\PPE$ instead of $X$ for cohomological considerations.

It is furthermore known, that the Picard group $\Pic(\PPE)$ is generated 
by the ruling $R=[\pi^*\OO_{\PP^1}(1)]$ and the hyperplane class
$H=[j^*\OO_{\PP^r}(1)]$ with intersection products
$$
H^d=f, \ \ H^{d-1}\cdot R=1, \ \ R^2=0.
$$
Hence, we will write a line bundle $\OO_\PPE(aH+bR)$ in the following form 
$$
\OO_\PPE(aH+bR) = \pi^*(\OO_{\PP^1}(b))(aH).
$$

\begin{theorem}[\cite{Sch}, Corollary 4.4]\label{Sch4.4}
Let $C$ be a curve with a complete base point free $g^1_k$ and 
let $\PPE$ be the projective bundle associated to the scroll $X$, swept out by the $g^1_k$.
\begin{enumerate}
\item [(a)]
$C\subset \PPE$ has a resolution $F_\bullet$ of type 
$$
0 \to \pi^*N_{k-2}(-kH) \to \pi^*N_{k-3}((-k+2)H) \to \dots \to \pi^*N_1(-2H) \to \OO_\PPE \to \OO_C \to 0
$$
with $\displaystyle N_i=\bigoplus_{j=1}^{\beta_i}\OO_{\PP^1}(a_j^{(i)})$ and $\displaystyle\beta_i=\frac{i(k-2-i)}{k-1}\binom{k}{i+1}$. 
\item [(b)] The complex $F_\bullet$ is self dual, i.e., 
$\cHom(F_\bullet, \OO_\PPE(-kH+(f-2)R))\cong F_\bullet$
\end{enumerate}
\end{theorem}

The resolution $F_\bullet$ above is called the \emph{relative canonical resolution}.
The degree of the bundles $N_i$ was computed in \cite{BH}.

\begin{proposition}[\cite{BH}, Proposition 2.9]\label{bundleDegs}
The degree of the bundle $N_i$ of rank $ \beta_i=\frac{k}{i+1}(k-2-i)\binom{k-2}{i-1}$ 
in the relative canonical resolution $F_\bullet$ is 
\vspace{-2mm}
$$
\deg(N_i)=\sum_{j=1}^{\beta_i}a_j^{(i)}=(g-k-1)(k-2-i)\binom{k-2}{i-1}.
$$
\end{proposition}

Since the rank and degree of the syzygy bundles $N_i$ over $\PP^1$ are known, the main object of investigation is the \emph{splitting type}. 

\begin{remark}\label{resolutionCE}
In \cite{CE}, Casnati and Ekedahl generalized the relative canonical resolution to finite Gorenstein covers $\pi:X\to Y$ of degree $k$. 
They define a relative resolution of $X\subset \PP(\cE_T)$, where $\cE_T$ is the 
the Tschirnhausen bundle on $Y$ defined by short exact sequence
$$
0\to \OO_Y \to \pi_*(\OO_X) \to \cE_T^{\vee} \to 0.
$$
Note that for a cover $C \overset{k:1}{\longrightarrow} \PP^1$, $\cE_T=\cE\otimes \OO_{\PP^1}(2)$, where $\cE$ is the bundle associate to $(C,g^1_k)$ as in Theorem \ref{Sch4.4}. 
The twists and hence the splitting types of the syzygy bundles in a resolution of $C\subset \PP(\cE_T)$ also differ from the ones in the relative canonical resolution of $C\subset \PPE$.
Indeed, following the proof of \cite[Step B, p. 445]{CE}, each twist in the $i$-th syzygy bundles in a resolution of $C\subset \PP(\cE_T)$ differs by exactly $2\cdot (i+1)$ from the ones given in our definition. 
Hence, we can deduce the degrees of the bundles  in this relative resolution of $C\subset \PP(\cE_T)$ from Proposition \ref{bundleDegs}. These degrees have also been computed directly in \cite{DP18}.   

\end{remark}

\begin{definition}
We say that a bundle on $\PP^1$ of the form $ N=\bigoplus_{j=1}^{\beta}\OO_{\PP^1}(a_j)$ is \emph{balanced} if
$$\displaystyle\max_{i,j} |a_j-a_i|\leq 1.$$ 
Equivalently, the bundle $N$ is balanced if $h^1(\PP^1,\cE\! nd(N))=0$. 
The relative canonical resolution is called balanced if 
all bundles $N_i$ occurring in the resolution are balanced.
\end{definition}

\begin{remark}
The locus of curves inside $\cH_{g,k}$  which have a balanced relative canonical resolution forms an open subset of $\cH_{g,k}$  which might be empty.
Hence, to show the generic balancedness for fixed values $(g,k)$ it is sufficient to examine a single balanced example. 
\end{remark}

\begin{remark}
The scroll associated to a general element in $\cH_{g,k}$ is always balanced by \cite{Bal} and \cite{H}. The sublocus inside $\cH_{g,k}$ parametrizing covers such that the associated scroll is unbalanced defines a divisor on $\cH_{g,k}$ precisely if $g$ is a multiple of $(k-1)$. This divisor is called the \emph{Maroni-divisor} (for more details on the Maroni-divisor see e.g. \cite{vdGK} and \cite{DP}). \\
On the other hand, knowing the splitting type of the syzygy bundles in the relative canonical resolution for generic elements in $\cH_{g,k}$ one can study the sublocus inside $\cH_{g,k}$ consisting set-theoretically of curves for which a certain syzygy bundle has non-generic splitting type. This yields interesting subvarieties which also turn out to be divisors in some cases (see \cite{DP18}). Similar to the Koszul-divisor on the moduli space $\cM_g$, the study of the divisors obtained from the relative canonical resolution sheds light on the global geometry of the Hurwitz space.
\end{remark}

\section{\texttt{Macaulay2}-package}

The \texttt{Macaulay2}-Package \texttt{RelativeCanonicalResolution.m2} (see \cite{M2} and \cite{BH4}) includes various useful functions to do experiments with $k$-gonal canonical curves and the relative canonical resolution of those curves.
We will briefly explain how functions in this package construct $g$-nodal $k$-gonal canonical curves of genus $g$. 

The main idea is that we start with a rational normalization of the desired curve and a degree $k$ map from the normalization to $\PP^1$. In the next step we pairwise glue $2g$ chosen points on the normalization. If $\cL$ is a line bundle of degree $k$ on a $g$-nodal curve $C$ with rational normalization $\nu:\PP\to C$, then $\cL$ is given as $\nu(\cL)\cong \OO_{\PP^1}(k)$ together with gluing data between the residue class fields 
$$\frac{a_i}{b_i}: \OO_{\PP^1}(k)\otimes \Bbbk(P_i)\to \OO_{\PP^1}(k)\otimes \Bbbk(Q_i).$$

Let $S=\Bbbk[s,t]$ be the coordinate ring of $\PP^1$. We start over by choosing two forms  $f,h\in S_k$ of degree $k$ and $g$ points $R_i=(R_i^{(0)}:R_i^{(1)})\in \PP^1$  such that for all $i=1,\dots,g$ the determinant
$$
\det \begin{pmatrix}
f & R_i^{(0)}\\
h & R_i^{(1)}
\end{pmatrix}
= l_i^{0}\cdot l_i^{(1)}\cdot r_i
$$
has at least two linear factors $l_i^{(0)}$ and $l_i^{(1)}$. Note that this step might be hard to perform over a field $\Bbbk$ of characteristic $0$ and we therefore
work over a finite field.
We compute $2g$ points $P_i=V(l_i^{(0)})$ and $Q_i=V(l_i^{(1)})$ as the vanishing loci of these linear forms. We want to define  multipliers $\{a_i,b_i\}_{i=1,\dots,g}$ such that 
$$
b_i \cdot f(P_i)=a_i \cdot h (Q_i) \text{ and } b_i \cdot h(P_i)=a_i \cdot f (Q_i) \text{ for } i=1,\dots,g.
$$
By construction, we can choose $\{a_i,b_i\}_{i=1,\dots,g}$ to be $b_i=1$ and $a_i=\frac{f(P_i)}{f(Q_i)}=\frac{h(P_i)}{h(Q_i)}$. 
If we define
$$
q_i:=\det \begin{pmatrix}
s & P_i^{(0)} \\
t & P_i^{(1)}
\end{pmatrix} 
\cdot 
\det \begin{pmatrix}
s & Q_i^{(0)} \\
t & Q_i^{(1)}
\end{pmatrix} \ \text{ for } i=1,\dots,g,
$$
then a basis of $H^0(C,\omega_C)$ can be identified with 
$$
\Bigg \{ 
s_j:=\prod_{ i=1,i\neq j}^g q_i \Bigg   \}_{j=1,\dots,g}.
$$
This basis  $\{s_j\}_{j=1,\dots, g}$ can furthermore be modified in such a way that the scroll defined by the line bundle of degree $k$ will have a "normalized" form, i.e., the $2\times (g-k+1)$ matrix defining the scroll will consist of blocks of the form
$$
\begin{pmatrix}
t_i & t_{i+2} \\
t_{i+1}& t_{i+3}
\end{pmatrix},
$$
where $T=\Bbbk[t_0,\dots,t_{g-1}]$ is the coordinate ring of $\PP^{g-1}$.

In the package \texttt{RelativeCanonicalResolution} we also provide a function which 
describes the generators of $C$ in terms of elements of the Cox ring of the scroll $\PPE$.

\begin{remark}
 There is an explicit identification 
 $$
 H^0(\PPE, \OO_\PPE(aH+bR))\cong H^0(\PP^1, (S_a\cE)(b)) \text{ for } a\ge 0,
 $$
 where $S_a\cE$ is the $a^{th}$ symmetric power of the vector bundle $\cE$ (see \cite[(1.3)]{Sch}). This gives a description of the coordinate ring 
 $$
 R_{\PPE} = \bigoplus_{a,b\in \ZZ} H^0(\PPE, \OO_\PPE(aH+bR))
 $$
 of $\PPE$ as the Cox ring $\Bbbk[v, w,\varphi_0,\dots,\varphi_{d-1}]$ equipped with bigrading $\deg v = \deg w = (1,0)$ and $\deg \varphi_i = (e_1-e_{i+1},1)$. 
\end{remark}

Finally the relative canonical resolution of $C\subset \PPE$ can be computed by successively picking syzygies in correct degrees.

\begin{example}
We compute a nodal $6$-gonal canonical curve of genus $9$.	
{
\begin{verbatim}
i1 : loadPackage("RelativeCanonicalResolution")
i2 : g=9; -- the genus
i3 : k=6; -- the degree of the pencil
i4 : n=10000; -- characteristic: next prime number after n 
i5 : Ican=canCurveWithFixedScroll(g,k,n); -- the canonical curve
i6 : (dim Ican,genus Ican, degree Ican)
o6 = (2, 9, 16)
i7 : betti(res(Ican,DegreeLimit=>1))
            0  1  2  3
o7 = total: 1 15 35 21
         0: 1  .  .  .
         1: . 21 64 70
\end{verbatim}	
}
Next we compute the ideal of $C$ inside the Cox ring of the scroll $\PPE$.
{
\begin{verbatim}
i8 : Jcan=curveOnScroll(Ican,g,k); -- the curve inside the scroll
i9 : RX=ring Jcan; -- the bigraded Cox ring of the scroll
       ZZ
o9 = -----[pp , pp , pp , pp , pp , v,w]
     10007   0    1    2    3    4
\end{verbatim}
}
We compute the relative canonical resolution:
{
\begin{verbatim}	
i10 : T=ring Ican; -- the canonical ring
i11 : H=basis({1,1},RX); -- a basis of H^0(PE, OO_PE(H))
i12 : phi=map(RX,T,H)
i13 : Ican==preimage_phi(Jcan)
o13 = true
i14 : lengthRes=2; -- a lengthlimit for the resolution on the scroll
\end{verbatim}
}
With respect to the total degree, the Betti table of the relative canonical resolution has the following form:
\newpage
{
\begin{verbatim}  
-- the relative canonical resolution:
i15 : betti(resX=resCurveOnScroll(Jcan,g,lengthRes)) 
             0 1  2 3 4 
o15 = total: 1 9 16 9 1
          0: 1 .  . . .
          1: . .  . . .
          2: . 6  2 . .
          3: . 3 12 3 .
          4: . .  2 6 .
          5: . .  . . 1
\end{verbatim}
}
The scroll cut out by the $g^1_6$ on $C$ has the following normalized determinantal representation:
{
\begin{verbatim}
i16 : X=preimage_phi(ideal 0_RX); -- the ideal of the scroll
i17 : repX=matrix{{t_0,t_2,t_4,t_6},{t_1,t_3,t_5,t_7}} 	
o17 = | t_0 t_2 t_4 t_6 |
      | t_1 t_3 t_5 t_7 |
i18 : minors(2,repX)==X
o18 = true      
\end{verbatim} 
}
\end{example}

\begin{remark}
By \texttt{o15}, we see that the second syzygy bundle $N_2$ is unbalanced in our example.
Although this single example does not show the generic unbalancedness for this case, one can show that this is indeed the generic form (see  \cite{BH2}).  	
\end{remark}

\section{Experiments and conjectures}

\subsection{Database of experiments}
Using our \texttt{Macaulay2}-package \texttt{RelativeCanonicalResolution.m2} we have computed the relative canonical resolution for various cases. For non-hyperelliptic, generic curves of genus $g\le 23$ with a pencil of degree $3\le k \le \min\{g-1,14\}$, all expected Betti tables are listed on the following webpage
 \vspace{0.5cm} \\
\resizebox{\textwidth}{!}{
 \url{https://www.math.uni-sb.de/ag/schreyer/images/data/computeralgebra/relcanres/html/index.html}
}
 \medskip \\
The webpage was set up with the help of Sascha Blug. 
All the experiments which led to Betti tables in \cite{M2exp} were performed over a finite field. If the examples for certain values $(g,k)$ yield a balanced relative canonical resolution, then by semi-continuity one can conclude that this is indeed the general behavior (even for complex algebraic curves). 

Since changing the characteristic for the unbalanced cases did not change the shape of the Betti tables, we believe that the Betti tables in \cite{M2exp} reflect the generic behavior. In general, we do not have a proof of this statement. However, the (un-)balancedness has been proved for the first bundle $N_1$ in some cases (see \cite{BH} and \cite{BP}). For several cases our examples lead to the conjecture, that certain higher syzygy bundles in the relative canonical resolution are unbalanced.  Most of these cases remain rather mysterious.

\subsection{Syzygy divisors on Hurwitz spaces}
Deopurkar and Patel used the relative canonical resolution to describe new effective divisors on the Hurwitz scheme $\cH_{g,k}$. If the degree $k$ divides $g-1$,
it is shown in \cite{BP} that the relative canonical resolution for a generic element in $\cH_{g,k}$ is totally balanced and hence, the locus $\mu_i$, corresponding set-theoretically to covers  in $\cH_{g,k}$ for which the $i$-th syzygy bundle $N_i$ is unbalanced, has expected codimension one. 
 In \cite{DP18} the authors give these syzygy divisor $\mu_1,\dots, \mu_{k-3}$ a scheme structure and compute their classes in a partial compactification of the Hurwitz scheme $\widetilde{\cH}_{g,k}$.
 In their main theorem, they represent the divisor classes $[\mu_i]$ in terms of certain tautological classes $\kappa, \zeta$ and $\delta$ (see \cite[\S 2]{DP18} for the precise definition of those classes).


\begin{theorem}[{\cite[Theorem 1.1]{DP18}}]\label{theoremDP18}
 Suppose $k$ divides $g-1$. 
 Let $i$ be an integer with $1\le i\le k-3$. The locus $\mu_i\subset \widetilde{\cH}_{g,k}$ is an effective divisor whose class in $\Pic_\QQ(\widetilde{\cH}_{g,k})$ is given by 
 $$
 [\mu_i] = A_i\cdot \bigg(6(gk-6g + k + 6)\cdot \zeta - k(k-12)\cdot \kappa - k^2\cdot \delta\bigg),
 $$
 where 
 $$
 A_i = \frac{(k-2)(k-3)}{6(i+1)(k-i-1)}\cdot \binom{k-4}{i-1}^2. 
 $$
\end{theorem}

Note that all the classes $[\mu_i]$ are proportional. The same phenomenon appears for classes of divisorial Brill--Noether loci in the moduli space $\overline{\cM}_g$. 
For the divisorial Brill--Noether classes it is known that these classes are supported on different sets and in  \cite{DP18} the authors conjecture that this also happens for the syzygy divisors on $\cH_{g,k}$. 

We come to a different conclusion.
Computing various examples of curves and their relative canonical resolution for $(g,k)\in \{(6,13),(7,15),(8,17),(6,19) \}$ over a field of small characteristic $p\le 500$ we found the following pattern which we conjecture to be true in general.
Note that computing random examples, the probability to end up in a certain codimension one locus is roughly $\frac{1}{p}$.


\begin{conjecture}
 Let $n,k$ be integers and $g-1 = n\cdot k$. Let $i$ be an integer with $1\le i\le k-3$. For a general element $(C,g^1_k)\in \mu_i \subset \widetilde{\cH}_{g,k}$, let $N_j$ be the $j$-th syzygy bundle in the relative canonical resolution of $C$ with $1\le j \le k-3$. Then all $N_j$ are unbalanced and the splitting type of $N_j$ is 
 $$
  N_j = \OO_{\PP^1}\big((n-1)(j+1)-1\big)^{\oplus \binom{k-4}{j-1}} \bigoplus \OO_{\PP^1}\big((n-1)(j+1)\big)^{\rk N_j - 2\cdot \binom{k-4}{j-1}} \bigoplus \OO_{\PP^1}\big((n-1)(j+1)+1\big)^{\oplus \binom{k-4}{j-1}}.
 $$
 In particular, all the effective divisors $\mu_i$ are supported on the same set. 
\end{conjecture}

\begin{remark}
One can easily check that the number $A_i$ in Theorem \ref{theoremDP18} is precisely
$$
A_i=\frac{1}{6k}\cdot \rk N_i \cdot \binom{k-4}{i-1}
$$
 The conjecture above predicts that the factor $\binom{k-4}{i-1}$ of $A_i$ also measures the unbalancedness of the bundle $N_i$. 
\end{remark}

\begin{remark}
If $(g-1)\neq n\cdot k$ then one can still consider the jumping loci set-theoretically defined as the subset of $\widetilde{\cH}_{g,k}$ consisting of covers such that the $i$-th syzygy bundle in the relative canonical resolution does not have generic splitting type. Similarly to the divisorial case one could ask if all those loci are supported on the same set. Experiments using our package \cite{BH4} show that there are several examples where these jumping loci have different support. 
\end{remark}

\subsection{Further conjectures}

We state several conjectures concerning the shape of relative canonical resolutions. This has partly also been discussed in \cite{BH}.

\begin{conjecture}\label{Conj-refinedBalancedness}
Let $C\subset \PP^{g-1}$ be a general canonical curve and let $k$ be a positive integer such that $\rho:=\rho(g,k,1)\geq 0$ and let $g^1_k$ be a general pencil in $W^1_k(C)$. Then for bundles $N_i=\bigoplus \OO_{\PP^1}(a_j^{(i)})$, $i=2,\dots, \big \lceil\frac{k-3}{2} \big \rceil$ there is the following bound 
$$
\max_{j,l} \big | a_j^{(i)}-a_l^{(i)} \big | \leq \min\{g-k-1,i+1\} \vspace{-1mm}.
$$ 
This bound is furthermore sharp in the following sense. Given two integers $k\geq 3$ and $2\leq i\leq  \lceil\frac{k-3}{2} \rceil$, there exists an integer $g$ such that the general canonical curve $C$ of genus $g$ has an $i$-th syzygy bundle $N_i$ in the relative canonical resolution, associated to a general pencil in $W^1_k(C)$, which satisfies $max_{j,l} | a_j^{(i)}-a_l^{(i)} | = \min\{g-k-1,i+1\} \vspace{-1mm}$.
In particular, if $g-k = 2$, the relative canonical resolution is balanced. 
\end{conjecture}

\begin{remark}
The above conjecture in the case $g-k=2$ says that the bundles in the relative canonical resolution are of the following form: 
$$
N_i=\OO_{\PP^1}^{\oplus i\cdot \binom{g-4}{i+1}}\bigoplus \OO_{\PP^1}(1)^{\oplus (g-4-i)\cdot \binom{g-4}{g-3-i}}.
$$
Note that the Betti numbers $i\cdot \binom{k-2}{i+1}$ appearing in the conjecture are the Betti numbers of a rational normal curve of degree $k-2$.
\end{remark}

We also verified the following Conjecture \ref{conj_balancedness_neg-rho} for $g\leq 23$.

\begin{conjecture}\label{conj_balancedness_neg-rho}
For a general cover $C\to \PP^1$ in $\cH_{g,k}$  with $\rho(g,k,1)\leq 0$, the bundle $N_1$ is balanced. 
\end{conjecture}

\paragraph*{Acknowledgement}

We would like to thank Sascha Blug for setting up the webpage \cite{M2exp} which displays all the experimental data. This work is a contribution to the Project 1.7 of the SFB-TRR 195 ''Symbolic Tools in Mathematics and their Application`` of the German Research Foundation (DFG). 


\end{document}